\documentclass[a4paper]{amsart}  

\usepackage{amssymb} \usepackage{amscd} \usepackage{times}

\def\zbb{\mathbb{Z}}  
  
  \def\phi{\varphi}
 \def\p1{{\mathbb{P}^1_\zbb}}

\begin{document}

\title{Some Uniforme Estimates for Scalar Curvature Type Equations.}
\author{Samy Skander Bahoura}

\address{6, rue Ferdinand Flocon, 75018 Paris, France. }
              
\email{samybahoura@yahoo.fr, samybahoura@gmail.com} 

\date{}

\maketitle

\begin{abstract}

We consider the prescribed scalar curvature equation on an open set $ \Omega  $ of  $ {\mathbb R}^n $, $ \Delta u=Vu^{(n+2)/(n-2)}+u^{n/(n-2)} $ with $ V\in C^{1,\alpha} $ $ (0 <\alpha \leq 1 ) $ and we prove the inequality $ \sup_K u \times \inf_{\Omega} u \leq c $ where $ K $ is a compact set of $ \Omega $.
 
In dimension 4, we have an idea on the supremum of the solution of the prescribed scalar curvature if we control the infimum. For this case we suppose the scalar curvature $ C^{1,\alpha} $, $ (0 < \alpha \leq 1) $.

\end{abstract}

\bigskip

\bigskip

\begin{center}  1. INTRODUCTION.
\end{center}

\bigskip

In our work, we denote $ \Delta=-\nabla^i \nabla_i $ the Laplace-Beltrami operator in dimension $ n\geq 2 $.

\bigskip

Here, we study some a priori estimates of type $ \sup \times \inf $ for prescribed scalar curvature equations in dimensions 4 and 5, also for perturbed scalar curvature equations in all dimension $ n \geq 3 $.

\bigskip

The $ \sup \times \inf $ inequality is caracteristic of those equations like the usual harnack inequalities for harmonic functions.

\bigskip

Note that, the prescribed scalar curvature equation was studied lot of. We can find,  see for example, [1], [2], [3], [4], [6], [7], [9], [10], [11], [12], [13] and [14], lot of results about uniform estimates in dimensions $ n=2 $ and $ n\geq 3 $.

\bigskip

In dimension 2, the corresponding equation is:

$$ \Delta u=Ve^u \qquad  (E_0) $$

Note that, Shafrir, see $ [14] $, have obtained an inequality of type $ \sup u +C \inf u <c $ with only $ L^{\infty} $ assumption on $ V $.

\bigskip

To obtain exactly the estimate $ \sup u + \inf u < c $, Brezis, Li and Shafrir gave a lipschitzian condition on $ V $, see [3].. Later, Chen and Lin have proved that if $ V $ is uniformly h\"olderian we can obtain a $ \sup + \inf $ inequality, see [7].

\bigskip

In dimension $ n \geq 3 $, the prescribed curvature equation on general manifold $ M $, is:

$$ \Delta u + R_g u =Vu^{(n+2)/(n-2)} \qquad (E'_0) $$

When $ M = {\mathbb S}_n $, Li, has proved  a priori estimates for the solutions of the previous equation. He use the notion of simple isolated points and some flatness conditions on $ V $, see [9] and [10].

\bigskip

If we suppose $ n=3,4 $, we can find in [12] and [13] uniforme estimates for the energy and a $ \sup \times \inf $ inequality. Note that, in [13], Li and Zhu have proved the compactness of the solutions of the Yamabe Problem by using the positive mass theorem.

\bigskip

In [2], we can see (on a bounded domain of $ {\mathbb R}^4 $) that we have an uniform estimate for the solutions of the equation $ (E'_0) $ ( $ n=4 $ and euclidian case) if we control the infimum of those functions, with only Lipschitzian assumption on the prescribed scalar curvature $ V $.

\bigskip

Here we extend  some result of [2] to equations with nonlinear terms or with minimal condition on the prescribed scalar curvature.

\bigskip

For the eulidian case, Chen and Lin gave some a priori estimates for general equations:

$$ \Delta u=Vu^{(n+2)/(n-2)}+g(u), \qquad (E''_0) $$

with some assumption on $ g $ and the Li-flatness conditions on $ V $, see [6].

\bigskip

Here, we give some a priori estimates with some minimal conditions on the precribed curvature. First, for perturbed scalar curvature equation, in all dimensions $ n \geq 3 $. Second, for prescribed scalar curvature equation  in dimensions $ 4 $ and $ 5 $.

\bigskip

Note that, we have no assumption on energy. In our work, we use the {\it blow-up} analysis, the {\it moving-plane} method. The moving-plane method was developped by Gidas-Ni-Nirenberg, see [8].

\bigskip

\bigskip

\begin{center}  2. MAIN RESULTS.
\end{center}

\bigskip

\bigskip

We consider the prescribed scalar curvature equation perturbed by a non linear term:

$$ \Delta u = Vu^{(n+2)/(n-2)}+u^{n/(n-2)} \,\,\, {\rm on} \,\,\, \Omega \subset {\mathbb R}^n \qquad (E_1).$$

Where $ V \in C^{1,\alpha} $, $ 0 < \alpha \leq 1 $, $ 0 < a \leq V(x) \leq b $ and $ ||V||_{C^{1,\alpha}} \leq A $.

\bigskip

We have,

\bigskip

{\it {\bf Theorem 1.} For all $ a, b, A, \alpha > 0 $ $ (0 < \alpha \leq 1 $) and all compact set $ K $ of $ \Omega $, there is a positive constant $ c=c(a, b, A, \alpha, K, \Omega,n) $ such that:

$$ \sup_K u \times \inf_{\Omega} u \leq c. $$}

If we suppose $ V \in C^1(\Omega) $ and $ V \geq a > 0 $, we have,

\bigskip

{\it {\bf Theorem 2.} For all $ a>0 $, $ V $ and all compact $ K $ of $ \Omega $, there is a positive constant $ c=c(a, V, K, \Omega, n) $ such that:

$$ \sup_K u \times \inf_{\Omega} u \leq c, $$

for all solution $ u $ of $ (E_1) $ relatively to $ V $.}

\bigskip

Now, we suppose $ n= 4 $, and we consider the following equation (prescribed scalar curvature equation):

$$ \Delta u=Vu^3 \,\,\, {\rm on} \,\, \Omega \subset {\mathbb R}^4 \qquad (E_2) $$

with $ 0 < a \leq V(x) \leq b $ and $ ||V||_{C^{1,\alpha}} \leq A $, $ 0 < \alpha \leq 1 $.

\bigskip

We have:

\bigskip

{\it {\bf Theorem 3.} For all $ a, b, m, A, \alpha >0 $, $ (0 < \alpha \leq 1) $ and all compact $ K  $ of $ \Omega $, there is a positive constant $ c=c(a,b,m,A,\alpha,K,\Omega) $ such that:

$$ \sup_K u \leq c \,\,\, {\rm if} \,\, \min_{\Omega} u \geq m. $$}

If we suppose $ n=4 $ and $ V \in C^1(\Omega) $ and $ V \geq a >0 $ on $ \Omega $, we have:

\bigskip

{\it {\bf Theorem 4.} For all $ a, m >0 $, $ V \in C^1(\Omega) $ and all compact $ K \in \Omega $, there is a positive constant $ c=c(a, m, V, K, \Omega) $ such that:

$$ \sup_K u \leq c \,\,\, {\rm if} \,\,\min_{\Omega} u \geq m, $$

for all $ u $ solution of $ (E_2) $ relatively to $ V $.}

\newpage

\begin{center}  3. PROOFS OF THE THEOREMS.
\end{center}

\bigskip

\underbar {\bf Proof of the Theorems 1 and 2.}

\bigskip

\bigskip

\underbar {\it Proof of the Theorem 1}

\bigskip

Without loss of generality, we suppose $ \Omega = B_1 $ the unit ball of $ {\mathbb R}^n $. We want to prove an a priori estimate around 0.

\bigskip

Let $ (u_i) $ and $ (V_i) $ be a sequences of functions on $ \Omega $ such that:

$$ \Delta u_i = V_i {u_i}^{(n+2)/(n-2)}+u_i^{n/(n-2)}, \,\, u_i>0,  $$

with $ 0 < a \leq V_i(x) \leq b $ and $ ||V_i||_{C^{1,\alpha}} \leq A $.

\bigskip

We argue by contradiction and we suppose that the $ \sup \times \inf $ is not bounded.

\bigskip

We have:

\bigskip

$ \forall \,\, c,R >0 \,\, \exists \,\, u_{c,R} $ solution of $ (E_1) $ such that:

$$ R^{n-2} \sup_{B(0,R)} u_{c,R} \times \inf_M u_{c,R} \geq c, \qquad (H) $$

\underbar {\bf Proposition :}{\it (blow-up analysis)} 

\smallskip

There is a sequence of points $ (y_i)_i $, $ y_i \to 0 $ and two sequences of positive real numbers $ (l_i)_i, (L_i)_i $, $ l_i \to 0 $, $ L_i \to +\infty $, such that if we set $ v_i(y)=\dfrac{u_i(y+y_i)}{u_i(y_i)} $, we have:

$$ 0 < v_i(y) \leq  \beta_i \leq 2^{(n-2)/2}, \,\, \beta_i \to 1. $$

$$  v_i(y)  \to \left ( \dfrac{1}{1+{|y|^2}} \right )^{(n-2)/2}, \,\, {\rm uniformly \,\, on \,\, all \,\, compact \,\, set \,\, of } \,\, {\mathbb R}^n . $$

$$ l_i^{(n-2)/2} u_i(y_i) \times \inf_{B_1} u_i \to +\infty,$$

\underbar {\bf Proof of the proposition:}

\bigskip

We use the hypothesis $ (H) $, we take two sequences $ R_i>0, R_i \to 0 $ and $ c_i \to +\infty $, such that,

$$ {R_i}^{(n-2)} \sup_{B(0,R_i)} u_i \times \inf_{B_1} u_i \geq c_i \to +\infty, $$

Let $ x_i \in  { B(x_0,R_i)} $ be a point such that $ \sup_{B(0,R_i)} u_i=u_i(x_i) $ and $ s_i(x)=(R_i-|x-x_i|)^{(n-2)/2} u_i(x), x\in B(x_i, R_i) $. Then, $ x_i \to 0 $.

\bigskip

We have:

$$ \max_{B(x_i,R_i)} s_i(x)=s_i(y_i) \geq s_i(x_i)={R_i}^{(n-2)/2} u_i(x_i)\geq \sqrt {c_i}  \to + \infty. $$ 

We set:

$$ l_i=R_i-|y_i-x_i|,\,\, \bar u_i(y)= u_i(y_i+y),\,\,  v_i(z)=\dfrac{u_i [y_i+\left ( z/[u_i(y_i)]^{2/(n-2)} \right )] } {u_i(y_i)}. $$

Clearly we have, $ y_i \to x_0 $. We also obtain:

$$ L_i= \dfrac{l_i}{(c_i)^{1/2(n-2)}} [u_i(y_i)]^{2/(n-2)}=\dfrac{[s_i(y_i)]^{2/(n-2)}}{c_i^{1/2(n-2)}}\geq \dfrac{c_i^{1/(n-2)}}{c_i^{1/2(n-2)}}=c_i^{1/2(n-2)}\to +\infty. $$

\bigskip

If $ |z|\leq L_i $, then $ y=[y_i+z/ [u_i(y_i)]^{2/(n-2)}] \in B(y_i,\delta_i l_i) $ with $ \delta_i=\dfrac{1}{(c_i)^{1/2(n-2)}} $ and $ |y-y_i| < R_i-|y_i-x_i| $, thus, $ |y- x_i| < R_i $ and, $ s_i(y)\leq s_i(y_i) $. We can write:

$$ u_i(y) (R_i-|y-y_i|)^{(n-2)/2} \leq u_i(y_i) (l_i)^{(n-2)/2}. $$

But, $ |y -y_i| \leq \delta_i l_i $, $ R_i >l_i$ and $ R_i-|y- y_i| \geq R_i-\delta_i l_i>l_i-\delta_i l_i=l_i(1-\delta_i) $. We obtain,

$$ 0 < v_i(z)=\dfrac{u_i(y)}{u_i(y_i)} \leq \left [ \dfrac{l_i}{l_i(1-\delta_i)} \right ]^{(n-2)/2}\leq 2^{(n-2)/2} . $$

We set, $ \beta_i=\left ( \dfrac{1}{1-\delta_i} \right )^{(n-2)/2} $, clearly, we have, $ \beta_i \to 1 $.

\bigskip

The function $ v_i $ satisfies:

$$ \Delta v_i= \tilde V_i {v_i}^{(n+2)/(n-2)}+\dfrac{v_i^{n/(n-2)}}{[u_i(y_i)]^{2/(n-2)}} $$

where, $ \tilde V_i(y)=V_i \left [ y+y/[u_i(y_i)]^{2/(n-2)} \right ] $. Without loss of generality, we can suppose that $ \tilde V_i \to V(0)=n(n-2) $.

\bigskip

We use the elliptic estimates, Ascoli and Ladyzenskaya theorems to have the uniform convergence of $ (v_i) $  to $ v $ on compact set of $ {\mathbb R}^n $. The function $ v $ satisfies: 

$$ \Delta v=n(n-2)v^{N-1}, \,\, v(0)=1,\,\, 0 \leq v\leq 1\leq 2^{(n-2)/2}, $$

By the maximum principle, we have $ v>0 $ on $ {\mathbb R}^n $. If we use Caffarelli-Gidas-Spruck result, ( see [5]), we obtain, $ v(y)=\left ( \dfrac{1}{1+{|y|^2}} \right )^{(n-2)/2} $. We have the same properties that in [2].

\bigskip

{\bf Remark}. When we use the convergence on compact sets of the sequence $ (v_i) $, we can take an increasing sequence of compact sets and we see that, we can obtain, a sequence $ (\epsilon_i) $ such that $ \epsilon_i \to 0 $ and after we choose $ (\tilde R_i) $ such that $ \tilde R_i \to + \infty $ and finaly:

$$ \tilde R_i^{n-2} ||v_i-v||_{B(0,\tilde R_i)} \leq \epsilon_i. $$

We can say that we are in the case of the step 1 of the theorem 1.2 of [6].

\bigskip

\underbar {\bf Fundamental Point:}{\it (a consequence of the blow-up)}

\bigskip

According to the work of Chen-Lin, see step 2 of the proof of the theorem 1.3 in [6], in the blow-up point, the prescribed scalar curvature $ V $ is such that:

$$ \lim_{i\to +\infty} |\nabla V_i(y_i)|= 0 \qquad (P_0) $$

\underbar {\bf Polar Coordinates} {\it (Moving-Plane method)}

\bigskip

Now, we must use the same method than in the Theorem 1 of [2]. We will use the moving-plane method.

\bigskip

We must prove the lemma 2 of [2].

\bigskip

We set $ t\in ]-\infty, -\log2 ] $ and $ \theta \in {\mathbb S}_{n-1}
$ :

$$ w_i(t,\theta)=e^{(n-2)t/2}u_i(y_i+e^t\theta), \,\, {\rm and} \,\,
\bar V_i(t,\theta)=V_i(y_i+e^t\theta). $$

We consider the following operator $ L = \partial_{tt}-\Delta_{\sigma}-\dfrac{(n-2)^2}{4} $, with $
\Delta_{\sigma} $ the Laplace-Baltrami operator on $ {\mathbb
  S}_{n-1} $.

\bigskip

The function $ w_i $ satisfies the following equation:

$$ -Lw_i=\bar V_i{w_i}^{N-1}+e^{t} \times {w_i}^{n/(n-2)}. $$

For $ \lambda \leq 0 $, we set :

\bigskip

$ t^{\lambda}=2\lambda-t  $ $ w_i^{\lambda}(t,\theta)=w_i(t^{\lambda},\theta) $ and $ \bar
V_i^{\lambda}(t,\theta)=\bar V_i(t^{\lambda},\theta) $

\bigskip

First, like in [2], we have the following lemma:

\bigskip

\underbar{\bf Lemma 1:} 

\bigskip

Let $ A_{\lambda} $ be the following property:

\bigskip

$ A_{\lambda}=\{\lambda\leq 0,\,\,\exists \,\,
 (t_{\lambda},\theta_{\lambda })\in
   ]\lambda,t_i]\times  {\mathbb S}_{n-1},\,\,   {\bar
   w_i}^{\lambda}(t_{\lambda},\theta_{\lambda})-{\bar
   w_i}(t_{\lambda},\theta_{\lambda} ) \geq 0\} . $

\bigskip

Then, there is $ \nu \leq 0$, such that for $\lambda \leq \nu $, $ A_{\lambda}  $ is not true.

\bigskip

\underbar {\bf Remark:} Here we choose $ t_i= \log \sqrt {l_i} $, where $ l_i $ is chooses as in the proposition.

\bigskip

Like in the proof of the Theorem 1 of [2], we want to prove the following lemma:

\underbar{\bf Lemma 2:} 

\bigskip

For $\lambda \leq 0 $ we have :

$$ {w_i}^{\lambda} -w_i \leq 0 \Rightarrow -L({w_i}^{\lambda} - w_i) \leq 0, $$

on $ ]\lambda,t_i]\times {\mathbb S}_{n-1} $.

\bigskip

Like in [2], we have:

\bigskip
  
\underbar{\bf A useful point: }

\bigskip

$ {\xi}_i= \sup \{\lambda \leq { \bar\lambda_i}=2+\log
  \eta_i, {w_i}^{\lambda} - w_i < 0 $, on $
  ]\lambda,t_i]\times {\mathbb S}_{n-1}
  \} $. The real $ \xi_i $ exists.
  
\bigskip

First:

$$ w_i(2\xi_i-t,\theta)=w_i[(\xi_i-t+\xi_i-\log\eta_i-2)+(\log
\eta_i+2)] , $$

\underbar {\bf Proof of the Lemma 2:}

\bigskip

In fact, for each $ i $ we have $ \lambda=
\xi_i \leq \log \eta_i+2 $, ($ \eta_i=[u_i(y_i)]^{(-2)/(n-2)}) $. 

\bigskip

Note that,

$$ w_i(2\xi_i-t,\theta)=w_i[(\xi_i-t+\xi_i-\log\eta_i-2)+(\log
\eta_i+2)] , $$

if we use the definition of $ w_i $ then for $ \xi_i \leq t $:

 $$ w_i(2\xi_i-t,\theta)=e^{[(n-2)(\xi_i-t+\xi_i-\log\eta_i-2)]/2}e^{n-2}v_i[\theta e^2e^{(\xi_i-t)+(\xi_i-\log\eta_i-2)}]
\leq 2^{(n-2)/2}e^{n-2}=\bar c. $$

We know that,

$$ -L( w_i^{\xi_i}-w_i)=[\bar V_i^{\xi_i
  }(w_i^{\xi_i})^{(n+2)/(n-2)}-\bar V_i
{w_i}^{(n+2)/(n-2)}]+[e^{t^{\xi_i}}(w_i^{\xi_i})^{n/(n-2) }-e^{t} {w_i}^{n/(n-2)} ],  $$ 

We denote by $ Z_1 $ and $ Z_2 $ the following terms:

$$ Z_1=(\bar V_i^{\xi_i }-\bar V_i)(w_i^{\xi_i })^{(n+2)/(n-2)}+\bar
V_i[(w_i^{\xi_i })^{(n+2)/(n-2)}-{w_i}^{(n+2)/(n-2)}],$$

and,

$$ Z_2 = e^{t^{\xi_i }}[(w_i^{\xi_i })^{n/(n-2)}-{w_i}^{n/(n-2)}]+ {w_i}^{n/(n-2)}(e^{t^{\xi_i }}-e^{t} ). $$

Like in the proof of the Theorem 1 of [2], we have:

$$ {w_i}^{\xi_i} \leq w_i \,\,\,{\rm and } \,\,\,
w_i^{\xi_i}(t,\theta)\leq \bar c \,\,\, {\rm for \, all } \,\,\,
(t,\theta)\in [\xi_i,-\log 2] \times {\mathbb S}_{n-1}  , $$

where, $ \bar c $ is a positive constant independant of $ i
$ and $ w_i^{\xi_i} $ for $ \xi_i \leq \log \eta_i+2 $.

\bigskip

\underbar {\it The $ (P_0) $ hypothesis:}

\bigskip

Now we use $ (P_0) $. We write:

$$ |\nabla V_i(y_i+e^t\theta)-\nabla V_i(y_i)|\leq A e^{\alpha t}, $$

Thus,

$$ |V_i(y_i+e^{t^{\xi_i}} \theta ) - V_i(y_i+e^t\theta)-<\nabla V_i(y_i)|\theta>(e^{t^{\xi_i}}-e^t)|\leq \dfrac{A}{1+\alpha} [e^{(1+\alpha)t^{\xi_i}}-e^{(1+\alpha )t}], $$

Then,

$$ |V_i^{\xi_i}-V_i|\leq |o(1)|(e^t-e^{t^{\xi_i}}), $$

Thus,

\bigskip

$ Z_1 \leq |o(1)| ({w_i^{\xi_i}})^{(n+2)/(n-2)}  (e^t-e^{ t^{\xi_i }}) \,\,\,
{\rm and } \,\,\,  Z_2 \leq {(w_i^{\xi_i})}^{n/(n-2)} \times   (e^{t^{\xi_i
    }}-e^{t} ) 
  $.

\bigskip

Then,

$$ -L(w_i^{\xi_i}-w_i) \leq (w_i^{\xi_i})^{n/(n-2)}[ (|o(1)|{w_i^{\xi_i}}^{2/(n-2)}-1)  (e^t-e^{ t^{\xi_i }})] \leq 0.
$$

The lemma is proved.

\bigskip

We set:

$$ \xi_i =\sup \{ \mu \leq \log \eta_i+2, w_i
  ^{\mu }(t,\theta)-w_i(t,\theta)  \leq 0, \forall \, (t,\theta) \in
  [\mu_i,t_i]\times {\mathbb S}_{n-1} \}, $$

with $ t_0 $ small enough.

\bigskip

Like in the proof of the Theorem 1 of [2], the maximum principle imply:

\smallskip

$$  \min_{\theta \in {\mathbb S}_{n-1}}w_i(t_i,\theta) \leq \max_{\theta
  \in {\mathbb S}_{n-1}} w_i(2\xi_i-t_i). $$

But,

$$ w_i(t_i,\theta)=e^{t_i} u_i(y_i+e^{t_i}\theta)\geq e^{t_i} \min u_i
\,\,{\rm and} \,\, w_i(2\xi_i-t_i)\leq \dfrac{ c_0 }{u_i(y_i)}, $$

thus,

$$ {l_i}^{(n-2)/2} u_i(y_i) \times \min u_i \leq c .$$

The proposition is contradicted.

\bigskip

\underbar {\it Proof of the Theorem 2.}

\bigskip

The proof of the Theorem 2 is similar to the proof of the Theorem 1. Only the "Fundamental point" change. We have:

\bigskip

According to the work of Chen-Lin, see step 2 of the proof of the theorem 1.1 in [6], in the blow-up point, the prescribed scalar curvature $ V $ is such that:

$$ \nabla V(0) = 0. $$

The function $ \nabla V $ is continuous on $ B_r(0) $ $ ( r $ small enough $ ) $, then it is uniformly continuous and we write (because $ y_i \to 0 $):

$$ |\nabla V(y_i+ y )-\nabla V(y_i) |\leq \epsilon, \,\,\, {\rm for } \,\,\, |y| \leq \delta <<r \,\, \forall \,\, i $$

Thus,

$$ |V^{\xi_i}-V|\leq o(1)(e^t-e^{t^{\xi_i}}), $$

We see that we have the same computations than in the section "Polar Coordinates" in the proof of the Theorem 1.

\bigskip

\underbar {\bf Proof of the Theorems 3 and 4.}

\bigskip

Here, only the section "Polar coordinates" change, the proposition of the first theorem stay true. First, we have:

\bigskip

\underbar {\bf Fundamental Point:}{\it (a consequence of the blow-up)}

According to the work of Chen-Lin, see step 2 of the proof of the theorem 1.3 in [6], in the blow-up point, the prescribed scalar curvature $ V $ is such that:

\bigskip

\underbar {\it Case 1: Theorem 3.}

$$ \lim_{i\to +\infty} |\nabla V_i(y_i)|= 0. $$

We write:

$$ |\nabla V_i(y_i+e^t\theta)-\nabla V_i(y_i)|\leq A e^{\alpha t}, $$

Thus,

$$ |V_i^{\xi_i}-V_i|\leq |o(1)|(e^t-e^{t^{\xi_i}}). $$

\underbar { \it Case 2: Theorem 4.}

$$ \nabla V(0) = 0. $$

The function $ \nabla V $ is continuous on $ B_r(0) $ $ ( r $ small enough $ ) $, then it is uniformly continuous and we write (because $ y_i \to 0 $):

$$ |\nabla V(y_i+ y )-\nabla V(y_i) |\leq \epsilon, \,\,\, {\rm for } \,\,\, |y| \leq \delta <<r \,\, \forall \,\, i $$

Thus,

$$ |V^{\xi_i}-V|\leq o(1)(e^t-e^{t^{\xi_i}}), $$

\underbar {\it Conclusion for Theorems 3 and 4.}

\bigskip

Finaly, we can note that we are in the case of the Theorem 2 of [2]. We have the same computations if we consider the following function:

$$ \bar w_i(t,\theta)= w_i(t,\theta)-\dfrac{m}{2} e^t. $$

We set, $ L= \partial_{tt}-\Delta_{\sigma }+1 $, where $ \Delta_{\sigma} $ is the Laplace-Beltrami operator on $ {\mathbb S}_3 $ and $ \bar V_i(t,\theta)=V_i(y_i+e^t\theta) $.

\bigskip

Like in [2], we want to prove the following lemma:

\bigskip

\underbar {Lemma.}

\smallskip

$$ \bar w_i^{\xi_i}-\bar w_i  \leq 0 \Rightarrow - L(\tilde
w_i^{\xi_i}-\tilde w_i)\leq 0 .$$

\underbar {Proof of the Lemma.}

$$ -L(\bar w_i^{\xi_i}-\bar w_i)=\bar V_i^{\xi_i}{(
  w_i^{\xi_i})}^3 -\bar V_i{ w_i}^3 .$$

Then,

$$ -L(\bar w_i^{\xi_i}-\bar w_i)=(\bar V_i^{\xi_i}-
 \bar V_i) {(w_i^{\xi_i})}^3+[
  {(w_i^{\xi_i})}^3 -{ w_i}^3]\bar V_i  .$$

For $ t\in [\xi_i,t_i] $ and $ \theta \in {\mathbb
  S}_3 $ :

$$ |\bar V_i^{\xi_i}(t,\theta)-\bar
V_i(t,\theta)|= |V_i(y_i+e^{2\xi_i-t}\theta )-V_i(y_i+e^t\theta)| \leq |o(1)|
(e^t-e^{2\xi_i-t}) . $$ 

The real $ t_i=\log \sqrt {l_i} \to -\infty $, where $ l_i $ is chooses as in the proposition of the theorem 1.

\bigskip

But, if $ \bar
  w_i^{\xi_i} -{ \bar w_i} \leq 0 $, we obtain:

$$  w_i^{\xi_i} - w_i \leq 
\dfrac{m}{2}(e^{2\xi_i-t}-e^t) <0  .$$

We use the fact that $ 0 < w_i^{\xi_i}< w_i $, we have:

$$ (w_i^{\xi_i})^3 -{ w_i}^3=(w_i^{\xi_i}
-w_i)[(w_i^{\xi_i})^2+w_i^{\xi_i} w_i+(w_i)^2] \leq 3 (w_i^{\xi_i}
-w_i) \times (w_i^{\xi_i})^2   .$$

Thus,we have for $ t\in [\xi_i,t_i] $ and $ \theta \in {\mathbb
  S}_3 $ :

$$  (w_i^{\xi_i})^3 -{ w_i}^3 \leq 3 \dfrac{m}{2}\, (w_i^{\xi_i})^2 (e^{2\xi_i-t}-e^t) . $$

We can write,

$$ -L(\bar w_i^{\xi_i}-\bar w_i)\leq (w_i^{\xi_i})^2 \,
(\dfrac{3 m}{2} \bar V_i-|o(1)|{w_i}^{\xi_i})\, (e^{2\xi_i-t }-e^t). \qquad (**)$$

We know that for $  t\leq \log
(l_i)-\log 2 +\log \eta_i $, we have,

$$ w_i(t,\theta)=e^t\times \dfrac{u_i\left ( y_i+\dfrac{e^t\theta}{u_i(y_i)}
  \right )}{u_i(y_i)} \leq 2 e^t . $$

We find,

$$ {w_i}^{\xi_i}(t,\theta) \leq 2e^2\sqrt {\dfrac{8}{a}},$$

because, $ \xi_i-\log \eta_i\leq 2+\dfrac{1}{2}\log \dfrac{8}{V(0)} $ and $ \xi_i\leq t\leq t_i $.

Finaly, $ (**) $ is negative and the Lemma is proved. 

\bigskip

Now, if we use the Hopf maximum principle, we obtain,

$$ \min_{\theta \in {\mathbb S}^3} \tilde w_i(t_i,\theta) \leq \max_{\theta
  \in {\mathbb S}^3}  \tilde w_i(2\xi_i-t_i,\theta) . $$

Which imply that,

$$ l_i u_i(y_i) \leq c. $$

It is a contradiction.

\bigskip

\bigskip

\underbar {\bf References:}

\bigskip

[1] T. Aubin. Some Nonlinear Problems in Riemannian Geometry. Springer-Verlag 1998.

\smallskip

[2] S.S Bahoura. Majorations du type $ \sup u \times \inf u \leq c $ pour l'\'equation de la courbure scalaire sur un ouvert de $ {\mathbb R}^n, n\geq 3 $. J. Math. Pures. Appl.(9) 83 2004 no, 9, 1109-1150.

\smallskip

[3] H. Brezis, Yy. Li Y-Y, I. Shafrir. A sup+inf inequality for some
nonlinear elliptic equations involving exponential
nonlinearities. J.Funct.Anal.115 (1993) 344-358.

\smallskip

[4] H.Brezis and F.Merle, Uniform estimates and blow-up bihavior for solutions of $ -\Delta u=Ve^u $ in two dimensions, Commun Partial Differential Equations 16 (1991), 1223-1253.

\smallskip

[5], L. Caffarelli, B. Gidas, J. Spruck. Asymptotic symmetry and local
behavior of semilinear elliptic equations with critical Sobolev
growth. Comm. Pure Appl. Math. 37 (1984) 369-402.

\smallskip

[6] C-C.Chen, C-S. Lin. Estimates of the conformal scalar curvature
equation via the method of moving planes. Comm. Pure
Appl. Math. L(1997) 0971-1017.

\smallskip

[7] A sharp sup+inf inequality for a nonlinear elliptic equation in ${\mathbb R}^2$.
Commun. Anal. Geom. 6, No.1, 1-19 (1998).

\smallskip

[8] B. Gidas, W-M. Ni, L. Nirenberg. Symmetry and Related Properties via the Maximum Principle. Commun. Math. Phys. 68, 209-243 (1979).

\smallskip

[9] YY. Li. Harnack Type Inequality: the Method of Moving Planes. Commun. Math. Phys. 200,421-444 (1999).

\smallskip

[10] YY. Li. Prescribing scalar curvature on $ {\mathbb S}_n $ and related
Problems. C.R. Acad. Sci. Paris 317 (1993) 159-164. Part
I: J. Differ. Equations 120 (1995) 319-410. Part II: Existence and
compactness. Comm. Pure Appl.Math.49 (1996) 541-597.

\smallskip

[11] YY. Li, I. Shafrir. Blow-up Analysis for Solutions of $ -\Delta u=Ve^u $ in Dimension Two. Indiana University Mathematics Journal. Vol. 43, No 4, (1994) 1255-1270.

\smallskip

[12] YY. Li, L. Zhang. A Harnack type inequality for the Yamabe equation in low dimensions.  Calc. Var. Partial Differential Equations  20  (2004),  no. 2, 133--151.

\smallskip

[13] YY.Li, M. Zhu. Yamabe Type Equations On Three Dimensional Riemannian Manifolds. Commun.Contem.Mathematics, vol 1. No.1 (1999) 1-50.

\smallskip

[14] I. Shafrir. A sup+inf inequality for the equation $ -\Delta u=Ve^u $. C. R. Acad.Sci. Paris S\'er. I Math. 315 (1992), no. 2, 159-164.

\end{document}